\documentclass[letterpaper, 10 pt, journal, onecolumn]{IEEEtran}
\IEEEoverridecommandlockouts
\usepackage[utf8]{inputenc}
\usepackage{cite}
\usepackage{amsmath,amssymb,amsfonts}
\usepackage{algorithmic}
\usepackage{graphicx}
\usepackage{textcomp}

\usepackage{amsmath,bm}
\usepackage{subfig}

\usepackage{graphicx}

\newtheorem{Theorem}{Theorem}
\newtheorem{Lemma}{Lemma}
\newtheorem{Remark}{Remark}

\newtheorem{Problem}{Problem}

\newtheorem{Assumption}{Assumption}

\usepackage{xcolor}

\newcommand{\m}[1]{\boldsymbol{#1}}
\newcommand{\mc}[1]{\mathcal{#1}}
\newcommand{\mb}[1]{\mathbb{#1}}

\begin{document}

\title{Finite-time bearing-based maneuver of acyclic leader-follower formations}
\author{Minh Hoang Trinh,~\IEEEmembership{Member,~IEEE,} and Hyo-Sung Ahn, \IEEEmembership{Senior Member,~IEEE}
\thanks{The work of M.H. Trinh is funded by the Hanoi University of Science and Technology (HUST) under project number T2020-SAHEP-007. The work of H.-S. Ahn was supported by the National Research Foundation of Korea (NRF) under the grant NRF- 2017R1A2B3007034}
\thanks{M. H. Trinh is with Department of Automatic Control, School of Electrical Engineering, Hanoi University of Science and Technology (HUST), Hanoi 11615, Vietnam. {E-mail: minh.trinhhoang@hust.edu.vn}}
\thanks{H.-S. Ahn is with School of Mechanical Engineering,  Gwangju Institute of Science and Technology (GIST), Gwangju 61005, Korea.
{E-mail: hyosung@gist.ac.kr}.}}

\maketitle
\thispagestyle{empty}
\begin{abstract}
This letter proposes two finite-time bearing-based  control laws for acyclic leader-follower formations. The leaders in formation move with a bounded continuous reference velocity and each follower controls its position with regard to three agents in the formation. The first control law uses only bearing vectors, and finite-time convergence is achieved by properly selecting two state-dependent control gains. The second control law requires both bearing vectors and communications between agents. Each agent simultaneously localizes and follows a virtual target. Finite-time convergence of the desired formation under both control laws is proved by mathematical induction and supported by numerical simulations.
\end{abstract}

\begin{IEEEkeywords}
bearing-only measurements, formation control, directed acyclic graph 
\end{IEEEkeywords}

\section{Introduction}
\label{sec1:introduction}
For many years, much research efforts have been put on understanding the mechanisms of collective behaviors displayed in nature and realizing them in large-scale systems such as robotic-, sensor-, traffic-, and electrical networks \cite{Anderson2008}. A notable application is formation control, where a team of autonomous agents (UAVs, UUVs, mobile robots, etc) is required to achieve and maintain a desired formation shape. Different solutions have been proposed to the problem based on various assumptions on the sensing/controlling/communication variables and topologies among the agents \cite{Oh2015}. 

In the bearing-based approach, the main sensing and controlling variables among the agents are the bearing vectors  (aka the directional information) or the subtended angles which can be obtained from low-cost cameras \cite{Eren2006using,Zhao2016aut,Zhao2018CSM}. Bearing-only formation control of a stationary formation with undirected \cite{Zhao2016tac,Tran2019finite} and special directed topologies  such as acyclic leader-follower \cite{Eren2012,Trinh2019TAC} or directed cycle \cite{Ko2020} have been studied in the literature. In flocking control, the agents simultaneously form a desired formation and agree on their velocities. To achieve flocking behavior, the agents sense the relative geometric variables and the relative velocity with regard to a few followers \cite{Moshtagh2006vision,Tran2018CCTA,Jing2019multiagent}. Formation tracking is more demanding since it requires the agents to achieve a target formation and follow a few  leaders, whose velocities can be time-varying. The ability to maintain a moving target formation shape is crucial for engineering applications such as search-and-rescue, truck platooning, or flight maneuvering. The authors in \cite{Zhao2015TCNS} proposed a control law for single and double-integrator agents using the relative positions and relative velocities. Bearing-only formation tracking with constant leaders' velocity has been studied for single integrators \cite{Zhao2019bearing,Zhao2020LCSS}, double integrators \cite{Zhao2019bearing,Trinh2021Aut}, or robotics agents \cite{Li2018NeuroComputing,Zhao2019bearing,Li2020adaptive,Huang2021TCNS}. However, bearing-only formation tracking with time-varying leader's velocity has not been studied in the literature. 

This work considers the bearing-only maneuver problem for directed acyclic leader-follower formations where the leaders' velocity is a bounded continuous function, thus filling a gap in the literature. The directed acyclic leader-follower structure can somehow describe the V-shape formation in immigrating birds, where each bird only sees and follows several individuals (immediate leaders) in its sight \cite{Bajec2009organized}. Thus, the interacting graph has a hierarchical structure and the formation is led by a small number of leaders \cite{Ballerini2008interaction,Nagy2010hierarchical}. To ease the analysis, it is assumed that during the formation maneuver, the positions of the immediate leaders are not collinear. Each follower, modeled by a single-integrator, controls its position based on information obtained from exactly three agents. First, a finite time bearing-only formation maneuver control law is proposed. Finite-time convergence of the target formation is achieved by appropriately choosing two control gains based on the measured bearing vectors. Note that this method has not been introduced in existing formation tracking laws in the literature \cite{Cao2011TAC,Vu2020ICCAS,Vu2020ICERA}. Second, in case the agents can communicate with each others, an estimation-based control strategy is proposed. Each agent determines a virtual target point based on the desired bearing vectors and the received position estimates of the immediate leaders. Simultaneously, the agent tracks its target point and sends this information to its followers. Although this control strategy requires more information than just the sensed bearing vectors, the followers have some freedom to choose their own trajectories to reach the target point. Thus, collision avoidance between agents or obstacle  can be included under this approach. Finally, simulations are given to support the analysis.

The remainder of this paper is organized as follows. Section \ref{sec2:problem_formulation} contains background and formulates the problem studied in this paper. Section \ref{sec:formation_tracking} proposes and studies the bearing-only formation tracking law. The target-point based formation tracking strategy is investigated in  Section \ref{sec:target_point}. Section \ref{sec:simulation} contains simulation results and Section \ref{sec:conclusion} concludes the letter. 

\section{Preliminaries and problem formulation}
\label{sec2:problem_formulation}
\subsection{Preliminaries}
\subsubsection{Notations}
The set of real numbers is denoted by $\mb{R}$. Given $x\in \mb{R}$, the signum function $\text{sign}(x)$ takes value 1, 0, -1 if $x>0$, $x=0$, and $x<0$, respectively. For $\alpha>0$, $\text{sig}(x)^\alpha = \text{sign}(x)|x|^\alpha$, where $|x|$ is the absolute value of $x$. Also, let $\m{x} \in \mb{R}^d$, one defines $|\m{x}|^{\alpha} = \sum_{i=1}^d |x_i|^\alpha$. The kernel and image of a matrix $\m{A} \in \mb{R}^{n\times n}$ are denoted by $\text{ker}(\m{A})$ and ${\text{im}}(\m{A})$. The vec operator is defined as vec$(\m{a}_1,\ldots,\m{a}_n)=[\m{a}_1^\top, \ldots, \m{a}_{n}^\top]^\top$. Let $\m{A} \in \mb{R}^{d \times d}$ be a symmetric matrix, one uses $\lambda_i(\m{A})$ to denote the $i$-th smallest eigenvalue of $\m{A}$.

Let $\m{x}$ be a nonzero vector in $\mb{R}^d$, defining the projection matrix $\m{P}_{\m{x}} = \m{I}_d - \frac{\m{x}\m{x}^\top}{\|\m{x}\|^2}$. For any vector $\m{y} \in \mb{R}^d$, $\m{P}_{\m{x}} \m{y}$ is the projection of $\m{y}$ onto the orthogonal space of $\text{im}(\m{x})$. $\m{P}_{\m{x}}$ is symmetric, positive semidefinite, $\text{ker}(\m{P}_{\m{x}}) = {\text{im}}(\m{x})$, and in addition to a zero eigenvalue, $\m{P}_{\m{x}}$ has $n-1$ eigenvalues 1. 
\subsubsection{Graph theory}
A directed graph $\mc{G}$ consists of a vertex set $\mc{V} = \{1, \ldots, n\}$ of $|\mc{V}| = n$ vertices and an edge set $\mc{E} \subset \mc{V}\times \mc{V}$ of $|\mc{E}| = m$ edges. The neighbor set of a vertex $i \in \mc{V}$ is defined as $\mc{N}_i = \{j \in \mc{V}|~(i, j) \in \mc{E}\}$. A simple path is a sequence of edges in $\mc{E}$ connecting  vertices in $\mc{V}$ so that there is no  repeated vertex (excepting for possibly the first and the last vertices) and edges. If the first and the last vertices of a path coincide, it is called a directed cycle. $\mc{G}$ is a directed acyclic graph if it does not contain any directed cycle. Indexing the edges so that $\mc{E}=\{e_1, \ldots, e_m\}$, the incidence matrix $\m{H} = [h_{ki}] \in \mb{R}^{m \times n}$ has $h_{ki} = -1$ if $i$ is the starting vertex of $e_k$, $h_{ki} = 1$ if $i$ is the end vertex of $e_k$, and $h_{ki}=0$, otherwise.
\subsection{Problem formulation}
Consider a formation of $n$ agents in the $d$-dimensional space ($d\ge 2$). Each agent in the formation has a local coordinate system $^i\Sigma$, and the axes of these local coordinate systems are aligned to each other. The position of agent $i$, written in a global coordinate system $^g\Sigma$, is denoted by $\m{p}_i \in \mb{R}^d$. 

A formation is characterized by $(\mc{G},\m{p})$, where $\mc{G}=(\mc{V},\mc{E})$ represents both the sensing and control interactions among agents and $\m{p} = \text{vec}(\m{p}_1,\ldots,\m{p}_n) \in \mb{R}^{dn}$ is a configuration of the formation. If $\m{p}_i\neq \m{p}_j$, the bearing vector between $i$ and $j$ is defined as $\m{g}_{ij} = \frac{\m{p}_j -\m{p}_i}{\|\m{p}_j-\m{p}_i\|} = \frac{\m{z}_{ij}}{\|\m{z}_{ij}\|}$, where $\m{z}_{ij}$ is the displacement between agents $i$ and $j$. If $(i,j) \in \mc{E}$, agent $i$ can sense $\m{g}_{ij}$ and control a desired bearing vector $\m{g}_{ij}^*$ with regard to agent $j$. The set of desired bearing vectors in the formation is denoted by ${\Gamma}=\{\m{g}_{ij}^*\}_{(i,j)\in\mc{E}}$. It is assumed that $\Gamma$ is feasible, i.e., there exists a target configuration $\m{p}^* = \text{vec}(\m{p}_1^*,\ldots, \m{p}_n^*) \in \mathbb{R}^{dn}$ satisfying all the bearing vectors in $\Gamma$.

Following the definition in \cite{Zhao2016tac}, $(\mc{G},\m{p})$ and $(\mc{G},\m{q})$ are bearing equivalent if and only if $\m{P}_{(\m{p}_i-\m{p}_j)}(\m{q}_i - \m{q}_j) = \m{0}, \forall (i,j)\in \mc{E}$. They are bearing congruent if and only if $\m{P}_{(\m{p}_i-\m{p}_j)}(\m{q}_i - \m{q}_j) = \m{0}, \forall i,j\in \mc{V}, i\ne j$. A formation $(\mc{G},\m{p})$ is globally bearing rigid if any formation $(\mc{G},\m{q})$  bearing equivalent to $(\mc{G},\m{p})$ is also bearing congruent to it.

The augmented bearing rigidity matrix is defined as $\tilde{\m{R}}_B = \text{diag}(\m{P}_{\m{g}_k})(\m{H}\otimes \m{I}_d)$. The formation $(\mc{G},\m{p})$ is infinitesimally bearing rigid if and only if the kernel of $\tilde{\m{R}}_B$ is only spanned by infinitesimal bearing rigid motions, i.e., \[\text{ker}(\tilde{\m{R}}_B)=\text{im}(\m{1}_n\otimes \m{I}_d, \m{p})=\text{im}(\m{1}_n\otimes \m{I}_d, \m{p}-\m{1}_n\otimes \bar{\m{p}}),\] where $\bar{\m{p}} = (\m{1}_n^\top \otimes \m{I}_d)\m{p}/n$ is the formation centroid. 

Suppose that in the formation, there are $l \ge 3$ leaders moving under the following equation
\begin{align} \label{eq:leaders}
\dot{\m{p}}_i = \m{v}_i,~i\in \mc{V}_l = \{1, \ldots, l\},
\end{align}
where $\m{v}_i$ is the velocity of the leader $i$. {Denote $\m{p}^l = \text{vec}(\m{p}_1,\ldots, \m{p}_l)$, and $\m{v}^l=\text{vec}(\m{v}_1, \ldots, \m{v}_l)$.} The remaining agents are followers, which are modeled by single-integrators
\begin{align} \label{eq:followers}
\dot{\m{p}}_i = \m{u}_i,~\forall i \in \mc{V}_f=\{ l+1, \ldots, n\},
\end{align}
where $\m{p}_i, \m{u}_i \in \mb{R}^d$ are respectively the position and the control input of agent $i$. 

The following assumptions are adopted in this paper:
\begin{Assumption} \label{eq:assumption-1}
The leaders are not in collinear positions and they move with the same bounded continuous reference velocity $\m{v}^*(t)$ satisfying $\|\m{v}^{*}(t)\| < \beta$, $\forall t \ge 0$. Furthermore, no collision happens between the agents. 
\end{Assumption}

\begin{Assumption} \label{eq:assumption-2}
The followers cannot sense the leaders' velocity but have information on the upper bound $\beta$. The bearing sensing and controlling graph $\mc{G}$ is a directed graph generated by the following procedure:
\begin{itemize}
\item Starting with $l$ vertices $1, \ldots, l$ where $l\ge 3$.
\item For each $l+1 \le i \le n$, inserting vertex $i$ together with $r \ge 3$ new edges $(i,j_k)$, where $j_k\in \{1, \ldots, i-1\}$.
\end{itemize}
\end{Assumption}
Obviously, any graph $\mc{G}$ satisfies Assumption~\ref{eq:assumption-2} is directed acyclic (see Fig.~\ref{fig:Sim1}(a) for an example).

\begin{Assumption} \label{eq:assumption-3}
The set of desired bearing vectors $\Gamma$ is feasible. The desired bearing vectors of each agent $i \in \mc{V}_f$, given by $\m{g}_{ij}^*, \forall j \in \mc{N}_i$, are not all parallel to each other.
\end{Assumption}

In this paper, the following problems will be studied.
\begin{Problem} \label{Problem}
Let Assumptions \ref{eq:assumption-1}--\ref{eq:assumption-3} hold. Design control laws for the followers modeled by \eqref{eq:followers} using only bearing vector measurements so that $\m{g}_{ij} \to \m{g}_{ij}^*, \forall (i,j) \in \mc{E}$ in finite time.
\end{Problem}

\begin{Problem} \label{Problem2}
Suppose that Assumptions \ref{eq:assumption-1}--\ref{eq:assumption-3} hold. Further, suppose that if $(i,j) \in \mc{E}$, agent $j$ sends its position estimate $\hat{\m{p}}_j$ to agent $i$ and the leaders have information on their true positions. Design control laws for the followers modeled by \eqref{eq:followers} so that $\m{g}_{ij} \to \m{g}_{ij}^*, \forall (i,j) \in \mc{E}$ in finite time.
\end{Problem}

\section{Finite-time bearing-only formation maneuver}
\label{sec:formation_tracking}
\subsection{Proposed control law}
The following bearing-only control law is proposed to solve Problem~\ref{Problem}:
\begin{align}
\m{u}_i &= -k_{i1} \sum_{j\in \mc{N}_i} \m{P}_{\m{g}_{ij}} \text{sig}(\m{P}_{\m{g}_{ij}}\m{g}_{ij}^*)^\alpha \nonumber\\
&\qquad\qquad\qquad - k_{i2} \beta \sum_{j\in \mc{N}_i} \m{P}_{\m{g}_{ij}} \text{sign}(\m{P}_{\m{g}_{ij}}\m{g}_{ij}^*), \label{eq:control-law-1}
\end{align}
where $i\in \mc{V}_f$, $\alpha \in (0,1)$, $\m{M}_i = \sum_{j\in \mc{N}_i} \m{P}_{\m{g}_{ij}}$, $k_{i1} = \lambda_1(\m{M}_i)^{-\frac{\alpha+1}{2}}$, $k_{i2} = \lambda_1(\m{M}_i)^{-\frac{1}{2}}$. In \eqref{eq:control-law-1}, the functions $\text{sig}(\cdot)$ and $\text{sign}(\cdot)$ are defined component-wise, $\m{P}_{\m{g}_{ij}} = \m{I}_d - \m{g}_{ij}\m{g}_{ij}^\top \in \mb{R}^{d\times d}$ is the projection matrix associating with $\m{g}_{ij}$. 

Some remarks on the proposed control law \eqref{eq:control-law-1} are given as follows. First, $\m{u}_i$ consists of two terms: the first term is for controlling the bearing-vectors of the agents to the desired ones; the second term is included to reject the disturbances resulting from leaders' motions which vary the bearing vectors. Second, as the leaders are not collinear, $\m{M}_i$ is symmetric positive definite \cite{Trinh2019TAC}. One has $\lambda_{1}(\m{M}_i) \leq \lambda_{n}(\m{M}_i) = \|\m{M}_i\| = \|\sum_{j\in\mc{N}_i}\m{P}_{\m{g}_{ij}}\| \leq \sum_{j\in\mc{N}_i}\|\m{P}_{\m{g}_{ij}}\| = |\mc{N}_i|$. Thus $k_{i1} \in (|\mc{N}_i|^{-\frac{\alpha+1}{2}}, \infty)$ and $k_{i2} \in \left(|\mc{N}_i|^{-\frac{1}{2}}, \infty \right)$. The number of neighbors for each follower $i$ is at least three, so that if its neighbors are in a non-collinear configuration and no collision happens between agents, $\lambda_1(\m{M}_i)>0$.

\subsection{Convergence Analysis}
Since $\mc{G}$ has a directed acyclic structure, the convergence analysis begins from the first follower. It will be shown that the first follower can achieve the desired bearing vectors and follow the leaders in finite time. Then, finite-time convergence of the overall desired formation will be established based on mathematical induction.
\subsubsection{The first follower}
Consider agent $i=l+1$ moving under the control law \eqref{eq:control-law-1}. Define the subformation $(\mc{K}_l,\m{p}^l)$, where $\mc{K}_l$ is the complete graph of $l$ vertices. 
{From Assumption~\ref{eq:assumption-1}, $\m{v}^l=\m{1}_l\otimes \m{v}^*$ belongs to the space of infinitesimal bearing rigid motions of $(\mc{K}_l,\m{p}^l)$. Thus, the relative bearing vector between the leaders are maintained during maneuver.}

The desired position of follower $i$ is determined by equations $\m{P}_{\m{g}_{ij}^*} (\m{p}_i^* - \m{p}_j) = \m{0}$, $\forall j \in \mc{N}_i$. By summing up these equations, it follows that $\sum_{j\in \mc{N}_i}\m{P}_{\m{g}_{ij}^*} (\m{p}_i^* - \m{p}_j) = \m{0}$ or equivalently $\sum_{j\in \mc{N}_i}\m{P}_{\m{g}_{ij}^*} \m{p}_i^* = \sum_{j\in \mc{N}_i}\m{P}_{\m{g}_{ij}^*} \m{p}_j $. According to the Assumption~\ref{eq:assumption-1}, the leaders are not colinear, which implies that the matrix $\sum_{j\in \mc{N}_i}\m{P}_{\m{g}_{ij}^*}$ is positive definite. It follows that $\m{p}_i^* = \big(\sum_{j\in \mc{N}_i} \m{P}_{\m{g}_{ij}^*} \big)^{-1} \big(\sum_{j\in \mc{N}_i} \m{P}_{\m{g}_{ij}^*} \m{p}_j \big)$ and thus, $\dot{\m{p}}_i^* = {\m{v}}_i^*  = \big(\sum_{j\in \mc{N}_i} \m{P}_{\m{g}_{ij}^*} \big)^{-1} \big(\sum_{j\in \mc{N}_i} \m{P}_{\m{g}_{ij}^*} \dot{\m{p}}_j \big) = \m{v}^*$. 

\begin{Lemma} \label{lem:firstfollower} Suppose that Assumptions~\ref{eq:assumption-1}--\ref{eq:assumption-3} hold. Under the control law~ \eqref{eq:control-law-1}, for $i=l+1$, $\m{p}_i \to \m{p}_i^*$ in finite time.
\end{Lemma}
\begin{IEEEproof}
Consider the Lyapunov function $V_i = \frac{1}{2}\|\m{p}_i - \m{p}_i^*\|^2$, which is positive definite and radially unbounded. As the control law \eqref{eq:control-law-1} is discontinuous, solution of \eqref{eq:followers} is understood in Filippov's sense \cite{Shevitz1994}. One has $\dot{V}_{i} \in^{a.e.} \dot{\tilde{V}}_i =  \bigcap_{\m{\xi}\in \nabla{V}_i(\m{p}_i)} \m{\xi}^\top K[\dot{\m{p}}_i]$. Then, $\nabla{V}_i(\m{p}_i) = \m{p}_i - \m{p}_i^*$, and
\begin{align*}
\dot{\tilde{V}}_i&= (\m{p}_i - \m{p}_i^*)^\top \Big(-k_{i1} \sum_{j\in \mc{N}_i} \m{P}_{\m{g}_{ij}} \text{sig}(\m{P}_{\m{g}_{ij}} \m{g}_{ij}^*)^\alpha  \\
&\qquad \qquad - k_{i2} \beta \sum_{j\in \mc{N}_i} \m{P}_{\m{g}_{ij}} K[\text{sign}](\m{P}_{\m{g}_{ij}} \m{g}_{ij}^*) - {\m{v}}^* \Big)  \\
&= (\m{p}_i - \m{p}_i^*)^\top \Big(-k_{i1} \sum_{j\in \mc{N}_i} \m{P}_{\m{g}_{ij}} \text{sig}\Big(\m{P}_{\m{g}_{ij}} \frac{\m{p}_j - \m{p}_i^*}{\|\m{z}_{ij}^*\|} \Big)^\alpha \nonumber\\
&\qquad\qquad  - k_{i2} \beta  \sum_{j\in \mc{N}_i} \m{P}_{\m{g}_{ij}} K[\text{sign}](\m{P}_{\m{g}_{ij}} \m{g}_{ij}^*) - {\m{v}}^* \Big).
\end{align*}
Because $\m{P}_{\m{g}_{ij}}(\m{p}_j - \m{p}_i^*)=\m{P}_{\m{g}_{ij}}(\m{p}_j - \m{p}_i + \m{p}_i - \m{p}_i^*) = \m{P}_{\m{g}_{ij}}(\m{p}_i - \m{p}_i^*),$
it follows that 
\begin{align}
\dot{\tilde{V}} &= -k_{i1} (\m{p}_i - \m{p}_i^*)^\top \sum_{j\in \mc{N}_i} \m{P}_{\m{g}_{ij}} \text{sig}\left(\frac{\m{P}_{\m{g}_{ij}}}{\|\m{z}_{ij}^*\|} (\m{p}_i - \m{p}_i^*) \right)^\alpha \nonumber\\
&\qquad - k_{i2} \beta  (\m{p}_i - \m{p}_i^*)^\top \sum_{j\in \mc{N}_i} \m{P}_{\m{g}_{ij}} K[\text{sign}](\m{P}_{\m{g}_{ij}} \m{g}_{ij}^*) \nonumber\\
&\qquad\qquad - (\m{p}_i - \m{p}_i^*)^\top{\m{v}}^*\nonumber\\
&\leq -k_{i1} (\m{p}_i - \m{p}_i^*)^\top \sum_{j\in \mc{N}_i} \m{P}_{\m{g}_{ij}} \text{sig}\Big(\frac{1}{\|\m{z}_{ij}^*\|}\m{P}_{\m{g}_{ij}} (\m{p}_i - \m{p}_i^*) \Big)^\alpha \nonumber\\
&\qquad - k_{i2} \beta  (\m{p}_i - \m{p}_i^*)^\top \sum_{j\in \mc{N}_i} \m{P}_{\m{g}_{ij}} K[\text{sign}](\m{P}_{\m{g}_{ij}} \m{g}_{ij}^*) \nonumber\\
&\qquad\qquad + \|\m{p}_i - \m{p}_i^*\|\|{\m{v}}^*\| \label{eq:dotV1}
\end{align}
From \cite{Shevitz1994} that for any $x \in \mb{R}$, $xK[\text{sign}](x) = \{|x|\}$, one has
\begin{align*}
(\m{p}_i &- \m{p}_i^*)^\top \sum_{j\in \mc{N}_i} \m{P}_{\m{g}_{ij}} \text{sig}\Big(\frac{1}{\|\m{z}_{ij}^*\|}\m{P}_{\m{g}_{ij}} (\m{p}_i - \m{p}_i^*) \Big)^\alpha \\
&=\sum_{j\in \mc{N}_i} \frac{1}{\|\m{z}_{ij}^*\|^\alpha} |\m{P}_{\m{g}_{ij}} (\m{p}_i - \m{p}_i^*) |^{\alpha+1} \nonumber\\
&= \sum_{j\in \mc{N}_i} \frac{1}{\|\m{z}_{ij}^*\|^\alpha} |(\m{p}_i - \m{p}_i^*)^\top \m{P}_{\m{g}_{ij}} (\m{p}_i - \m{p}_i^*) |^{\frac{\alpha+1}{2}} \\
& \ge \frac{1}{\max_{j\in \mc{N}_i} (\|\m{z}_{ij}^*\|)^\alpha} \left|(\m{p}_i - \m{p}_i^*)^\top \m{M}_i (\m{p}_i - \m{p}_i^*) \right|^{\frac{\alpha+1}{2}} \\
&\ge \frac{1}{\max_{j\in \mc{N}_i} (\|\m{z}_{ij}^*\|)^\alpha} \big(\lambda_1 (\m{M}_i) \big)^{\frac{\alpha+1}{2}} \|\m{p}_i - \m{p}_i^*\|^{\alpha+1},
\end{align*}
where the first inequality follows from \cite{Wang2010finite}[Lemma 2], and 
\begin{align*}
(\m{p}_i &- \m{p}_i^*)^\top \sum_{j\in \mc{N}_i} \m{P}_{\m{g}_{ij}} K[\text{sign}](\m{P}_{\m{g}_{ij}} \m{g}_{ij}^*) \nonumber\\
&= \sum_{j\in \mc{N}_i} (\m{p}_i - \m{p}_i^*)^\top \m{P}_{\m{g}_{ij}} K[\text{sign}]\Big(\m{P}_{\m{g}_{ij}}(\m{p}_i - \m{p}_i^*)\Big) \\
&= \sum_{j\in \mc{N}_i} \|\m{P}_{\m{g}_{ij}}(\m{p}_i - \m{p}_i^*)\|_1\ge \left|\left| \sum_{j\in \mc{N}_i}\m{P}_{\m{g}_{ij}}(\m{p}_i - \m{p}_i^*) \right|\right|_1\\
&\ge \left|\left| \sum_{j\in \mc{N}_i}\m{P}_{\m{g}_{ij}}(\m{p}_i - \m{p}_i^*) \right|\right|\ge {\lambda_1(\m{M}_i)}^{\frac{1}{2}} \|\m{p}_i - \m{p}_i^*\|.
\end{align*}
Thus, one has
\begin{align}
\dot{V}_i 
&\le - \eta V_i^{\frac{\alpha+1}{2}} - \big(\beta-\|{\m{v}}^*\| \big) \|\m{p}_i - \m{p}_i^*\|,
\end{align}
where $\eta = \frac{2^{\frac{\alpha+1}{2}}}{\max_{j\in \mc{N}_i} (\|\m{z}_{ij}^*\|)^\alpha}$. Based on \cite{Bhat2000}, $\m{p}_i \to \m{p}_i^*$ in a finite time upper bounded by $T_{l+1} = \frac{2V_{l+1}(0)^{\frac{1-\alpha}{2}}}{\eta(1-\alpha)}$.
\end{IEEEproof}

\subsubsection{The n-agent formation}
For each agent $i\ge l+1$, the velocities of agents $j \in \mc{N}_i \subset \{1,\ldots, i-1\}$ are considered as external inputs to the dynamics of agent $i$. The convergence of the overall formation is given in the following theorem.

\begin{Theorem} \label{thm:1} Suppose that Assumptions \ref{eq:assumption-1}--\ref{eq:assumption-3} hold. Further, suppose that each agent and its neighbors are always non-collinear. Under control law \eqref{eq:control-law-1}, $\m{p}_i(t) \to \m{p}_i^*$ in finite time upper bounded by $T_i = T_{i-1} + \frac{2{V}_i(T_{i-1})^{\frac{1-\alpha}{2}}}{\eta(1-\alpha)}$, $i=l+1, \ldots, n$ and $T_{l}=0$, and $\m{u}_i$ is bounded $\forall i=l+1, \ldots, n$.
\end{Theorem}

\begin{IEEEproof} We show this theorem by mathematical induction on $i$. From Lemma \ref{lem:firstfollower}, the claim is true for $i=l+1$. Next, suppose that the claim holds until $i~(l+1 \leq i \leq n)$, then $\forall j\in \mc{N}_i$, during $[0,T_{j-1}]$, $\m{u}_j$ is bounded and $\m{p}_j = \m{p}_j^*$ for $t\ge T_{j}$, where $T_k = T_{k-1} + \frac{2{V}_i(T_{k-1})^{\frac{1-\alpha}{2}}}{\eta(1-\alpha)}$, $\forall k=l+1, \ldots, i-1$ and $T_{l}=0$. Since agent $i$ and its neighbor agents are not collinear and no collisions happen, $\m{M}_i=\sum_{j \in \mc{N}_i} \m{P}_{\m{g}_{ij}}$ is always positive definite and $k_{i1}, k_{i2}<\infty$ in $t\in [0,T_{i-1}]$. It follows that $\dot{\m{p}}_i(t)$ is bounded and thus $\|\m{p}_i(T_{i-1}) - \m{p}_i^*(T_{i-1})\|<\infty$.

Next, for $t\ge T_{i-1}$, consider the Lyapunov function $V_i = \frac{1}{2}\|\m{p}_i - \m{p}_i^*\|^2$.  Similar to Lemma~\ref{lem:firstfollower}, one has
\begin{align} \label{eq:boundedness}
\dot{V}_i &\le - \frac{\|\m{p}_i - \m{p}_i^*\|^{\alpha+1}}{\max_{j \in \mc{N}_i} \|\m{z}_{ij}^*\|^\alpha} - \|\m{p}_i -\m{p}_i^*\| (\beta-\|\m{v}^*\|).
\end{align}
This means $\m{p}_i(t) \to \m{p}_i^*$ in a finite time upper bounded by $\sum_{k=l+1}^{i}T_k$, where $T_i = \frac{2V_i(T_{i-1})^{\frac{1-\alpha}{2}}}{\eta(1-\alpha)}$.

Therefore, the claim is also true for $i=n$ by mathematical induction.
\end{IEEEproof}

\begin{Remark} After the target formation was achieved, the leaders can rescale the formation by adopting the reference velocity $\m{v}^l(t) \in \text{im}(\m{1}_l\otimes \m{I}_d, \m{p}^l(t))$, $\m{v}^{l}(t) = \m{v}^{l}_{\text{trans}}(t) + k_{\text{scale}}(t) \m{v}^{l}_{\text{scale}}(t)$, where $\m{v}^{l}_{\text{trans}}(t) = \m{1}_l \otimes \m{v}^*$ and $\m{v}^{l}_{\text{scale}}(t) = \frac{\m{p}^l - \bar{\m{p}}^l\otimes \m{1}_l}{\|\m{p}^l - \bar{\m{p}}^l\otimes \m{1}_l\|}$, $\bar{\m{p}}^l=(\m{1}_l\otimes\m{I}_d)\m{p}^l/l$, are the translational and scaling motions of the formation $(\mc{K}_l,\m{p}^l)$, respectively. Suppose that $\m{v}^{l}_{\text{trans}}$ is bounded and 
$k_{\text{scale}}(t)= \xi >0$, $t \in [t_1,t_2]$, $0\le t_1 <t_2<\infty,$ and $k_{\text{scale}}(t)=0$, otherwise. This assumption describes that the formation mainly moves forward and the rescaling process is occasionally performed during the maneuver, e.g., before and after traversing a narrow alley. For $t\ge t_2$, convergence of the desired formation can be shown as in Thm.~\ref{thm:1}.
\end{Remark}

\begin{Remark} Since collision between agents and collinearity of $\{\m{p}_j\}_{j\in\mc{N}_i}, \forall i=l+1,\ldots, n$ are excluded, one cannot conclude in Thm.~\ref{thm:1} about globally asymptotic stability of the target formation. Also, the mathematical induction cannot be established without finite time convergence of each follower since boundedness of $\|\m{p}_i - \m{p}_i^*\|,\forall t$ will not be guaranteed.
\end{Remark}
\section{Finite-time formation maneuver via target point localization}
\label{sec:target_point}
\subsection{Proposed control law}
In this section, we propose an control strategy to solve Problem 2. Instead of directly controlling the position based on the bearing errors, each follower estimates its desired position with regard to the leaders and track that point. The following finite-time control law is proposed:
\begin{align} 
\dot{\hat{\m{p}}}_i & = \sum_{j\in \mc{N}_i} \m{P}_{\m{g}_{ij}^*}\text{sig}\big(\m{P}_{\m{g}_{ij}^*}(\hat{\m{p}}_j -  \hat{\m{p}}_i)\big)^\alpha  +\gamma_i \sum_{j\in \mc{N}_i} \m{P}_{\m{g}_{ij}^*} \text{sign}\big(\m{P}_{\m{g}_{ij}^*}(\hat{\m{p}}_j - \hat{\m{p}}_i) \big),\label{eq:est1}\\
\dot{\m{p}}_i &= - \text{sig}\big(\m{p}_i - \hat{\m{p}}_i \big)^\alpha - \beta \text{sign}\big(\m{p}_i - \hat{\m{p}}_i \big), \label{eq:est2}
\end{align}
where $\alpha \in (0, 1)$ and $\gamma_i \ge \beta \left|\left| \sum_{j \in \mc{N}_i} \m{P}_{\m{g}_{ij}^*} \right|\right|^{-1}$. Due to the directed acyclic structure of the graph $\mc{G}$, agent $i$ can also calculate the target point directly from the information by $\m{p}_i^* = (\sum_{j\in \mc{N}_i} \m{P}_{\m{g}_{ij}^*})^{-1} \sum_{j\in \mc{N}_i} \m{P}_{\m{g}_{ij}^*}\m{p}_j^*$, and track that point under the control law \eqref{eq:est2}. In case the beacons are stationary, a fixed-time network localization law was proposed in \cite{Trinh2020ACC}.
\subsection{Convergence analysis}
Consider the first follower $i=l+1$. The following lemma will be proved.
\begin{Lemma} \label{lem:est} Let the assumptions of Problem~\ref{Problem2} hold, for  $i=l+1$, ${\hat{\m{p}}}_i(t) \to {{\m{p}}}_i^*$ and ${\m{p}}_i(t) \to {\m{p}}_i^*$ in finite time.
\end{Lemma}
\begin{IEEEproof}
First, consider the estimation dynamics \eqref{eq:est1}. Using the Lyapunov function $V = \frac{1}{2}\|\hat{\m{p}}_i - \m{p}^*_i\|^2$ and keeping in mind that $\hat{\m{p}}_j = \m{p}_j^*, \forall j \in \mc{N}_i$, one has
\begin{align}
\dot{V} & = (\hat{\m{p}}_i - \m{p}^*_i)^\top \Big(\m{P}_{\m{g}_{ij}^*}\sum_{j\in \mc{N}_i} \text{sig}\big(\m{P}_{\m{g}_{ij}^*}({\m{p}}_j^* - {\m{p}}_i^* + {\m{p}}_i^*- \hat{\m{p}}_i)\big)^\alpha  + \gamma_i \m{P}_{\m{g}_{ij}^*}\sum_{j\in \mc{N}_i} \text{sign}\big(\m{P}_{\m{g}_{ij}^*}({\m{p}}_j^* - {\m{p}}_i^* + {\m{p}}_i^* - \hat{\m{p}}_i)\big) - \m{v}^* \Big) \nonumber \\
&= - \sum_{j \in \mc{N}_i} | \m{P}_{\m{g}_{ij}^*}({\m{p}}_i^* - \hat{\m{p}}_i)|^{\alpha+1} - \gamma_i \sum_{j \in \mc{N}_i} \| \m{P}_{\m{g}_{ij}^*}({\m{p}}_i^* - \hat{\m{p}}_i)\|_1 - (\hat{\m{p}}_i - \m{p}^*_i)^\top \m{v}^* \nonumber\\
&\le - \left| ({\m{p}}_i^* - \hat{\m{p}}_i)^\top \sum_{j \in \mc{N}_i} \m{P}_{\m{g}_{ij}^*}({\m{p}}_i^* - \hat{\m{p}}_i) \right|^{\frac{\alpha+1}{2}} - \gamma_i \left|\left| \sum_{j \in \mc{N}_i} \m{P}_{\m{g}_{ij}^*}({\m{p}}_i^* - \hat{\m{p}}_i) \right|\right|_1  + \|\m{v}^*\| \|{\m{p}}_i^* - \hat{\m{p}}_i\| \nonumber\\
&\le - \eta_i | ({\m{p}}_i^* - \hat{\m{p}}_i)^\top ({\m{p}}_i^* - \hat{\m{p}}_i) |^{\frac{\alpha+1}{2}} - \gamma_i \left|\left| \sum_{j \in \mc{N}_i} \m{P}_{\m{g}_{ij}^*} \right|\right| \|{\m{p}}_i^* - \hat{\m{p}}_i\| + \|\m{v}^*\| \|{\m{p}}_i^* - \hat{\m{p}}_i\|\nonumber\\
&\le - {\chi_i} V^{\frac{\alpha+1}{2}} - \Big(\gamma_i\| \sum_{j \in \mc{N}_i} \m{P}_{\m{g}_{ij}^*} \| - \|\m{v}^*\| \Big) \|{\m{p}}_i^* - \hat{\m{p}}_i\| \label{eq:est-dotV3}
\end{align}
where $\eta_i = \lambda_1\big(\sum_{j \in \mc{N}_i} \m{P}_{\m{g}_{ij}^*}\big)^{\frac{\alpha+1}{2}}$ and $\chi_i = {\eta_i}2^{\frac{\alpha+1}{2}}$. It follows that ${\hat{\m{p}}}_i(t) \to {{\m{p}}}_i^*$ in finite time upper bounded by $T_i$, and $\hat{\m{v}}_i \to {\m{v}}_i^*$ in finite time. 

Next, consider the position tracking control law~\eqref{eq:est2} with the Lyapunov function $W = \frac{1}{2}\|\m{p}_i - {\hat{\m{p}}}_i\|^2$. For $t\leq  T_i$, 
\begin{align}
\dot{W} &= -(\m{p}_i - {\hat{\m{p}}}_i)^\top \text{sig}\big(\m{p}_i - \hat{\m{p}}_i \big)^\alpha -  \beta (\m{p}_i - {\hat{\m{p}}}_i)^\top\text{sign}\big(\m{p}_i - \hat{\m{p}}_i \big) - (\m{p}_i - {\hat{\m{p}}}_i)^\top {\hat{\m{v}}}_i \nonumber\\
&\leq -|\m{p}_i - \hat{\m{p}}_i|^{\alpha+1} - \beta \|\m{p}_i - {\hat{\m{p}}}_i\|_1 + \|\m{p}_i - {\hat{\m{p}}}_i\| \|{\hat{\m{v}}}_i\|, \label{eq:13}
\end{align}
which shows that $\|\m{p}_i - {\hat{\m{p}}}_i\|$ is globally ultimately bounded. Together with the boundedness of $\|{\hat{\m{p}}}_i - \m{p}_i^*\|$, it follows from the  triangle inequality that $\|\m{p}_i -\m{p}_i^*\|\le \|\m{p}_i - {\hat{\m{p}}}_i\| + \|{\hat{\m{p}}}_i - \m{p}_i^*\|$, i.e., $\|\m{p}_i -\m{p}_i^*\|$ is also bounded. 

Now, for $t\ge T_i$, $\hat{\m{p}}_i = \m{p}^*_i$ and $\|\hat{\m{v}}_i\| = \|{\m{v}}^*\|\leq \beta$, consider the function $W_1 = \frac{1}{2}\|\m{p}_i - {\hat{\m{p}}}_i^*\|^2$. Similar to \eqref{eq:13}, one finds that $\dot{W}_1 \leq -|\m{p}_i - \hat{\m{p}}_i|^{\alpha+1} - \|\m{p}_i - {\hat{\m{p}}}_i\| (\beta- \|{\hat{\m{v}}}_i\|) \leq -\zeta_i W_1^{\frac{\alpha+1}{2}}$, where $\zeta_i = 2^{\frac{\alpha+1}{2}}> 0$ is a positive constant. Thus, $\m{p}_i \to \m{p}_i^*$ in finite time.
\end{IEEEproof}

\begin{Theorem} \label{thm:est} Under the control laws \eqref{eq:est1}--\eqref{eq:est2}, the desired moving formation is achieved in finite time.
\end{Theorem}
\begin{IEEEproof}
The proof follows from Lemma~\ref{lem:est} and mathematical induction on $i=l+1$ to $n$.
\end{IEEEproof}

\begin{Remark}
Observe that \eqref{eq:est2} steers the agent to the virtual target along a curve. Let $d=2$ and assume that the agents are equipped with proximity sensors which can sense the distance to an obstacle located at $\m{p}_o$ within a small range $d_{max}$. To avoid collision with an obstacle, for $0 < d < d_{\max}$, \eqref{eq:est2} can be modified as follows: 
\begin{align} \label{eq:coav}
{\m{u}}_i &= - (1-\zeta) \Big(\text{sig}\big(\m{p}_i - \hat{\m{p}}_i \big)^\alpha + \beta \text{sign}\big(\m{p}_i - \hat{\m{p}}_i \big) \Big) + \zeta \big(k(\m{p}_i - {\m{p}}_o) + \m{g}_{io}^\perp\big),
\end{align}
where $\zeta = 1$ if $\|\m{p}_i-\m{p}_0\|< d$ and $\zeta = 0$ if $\|\m{p}_i-\m{p}_o\|\ge d$, $\m{g}_{io}^\perp = \m{J} \frac{\m{p}_o - \m{p}_i}{\|\m{p}_o - \m{p}_i\|}$, $\m{J} = \begin{bmatrix}
0 & 1\\ -1 & 0
\end{bmatrix}$, $k>1$ is a sufficiently large control gain.
\end{Remark}
\begin{IEEEproof}
Denote ${C} = \{\m{q} \in \mb{R}^2: \|\m{q} - \m{p}_o\| < d\}$ and $\partial C = \{\m{q} \in \mb{R}^2: \|\m{q} - \m{p}_o\| = d\}$. Consider the following scenarios:
\begin{itemize}
\item Case 1: The solution of $\dot{\m{p}}_i = - \Big(\text{sig}\big(\m{p}_i - \hat{\m{p}}_i \big)^\alpha + \beta \text{sign}\big(\m{p}_i - \hat{\m{p}}_i \big) \Big), {\m{p}}_i(t_0) = {\m{p}}_{i0}$, which is called $S$, does not intersect $C$. Then, $\zeta = 0, \forall t$ and the agent moves along the curve to the target under the control law: ${\m{u}}_i = - \Big(\text{sig}\big(\m{p}_i - \hat{\m{p}}_i \big)^\alpha + \beta \text{sign}\big(\m{p}_i - \hat{\m{p}}_i \big) \Big)$.
\begin{figure*}[ht]
\centering
\includegraphics[height=3.5cm]{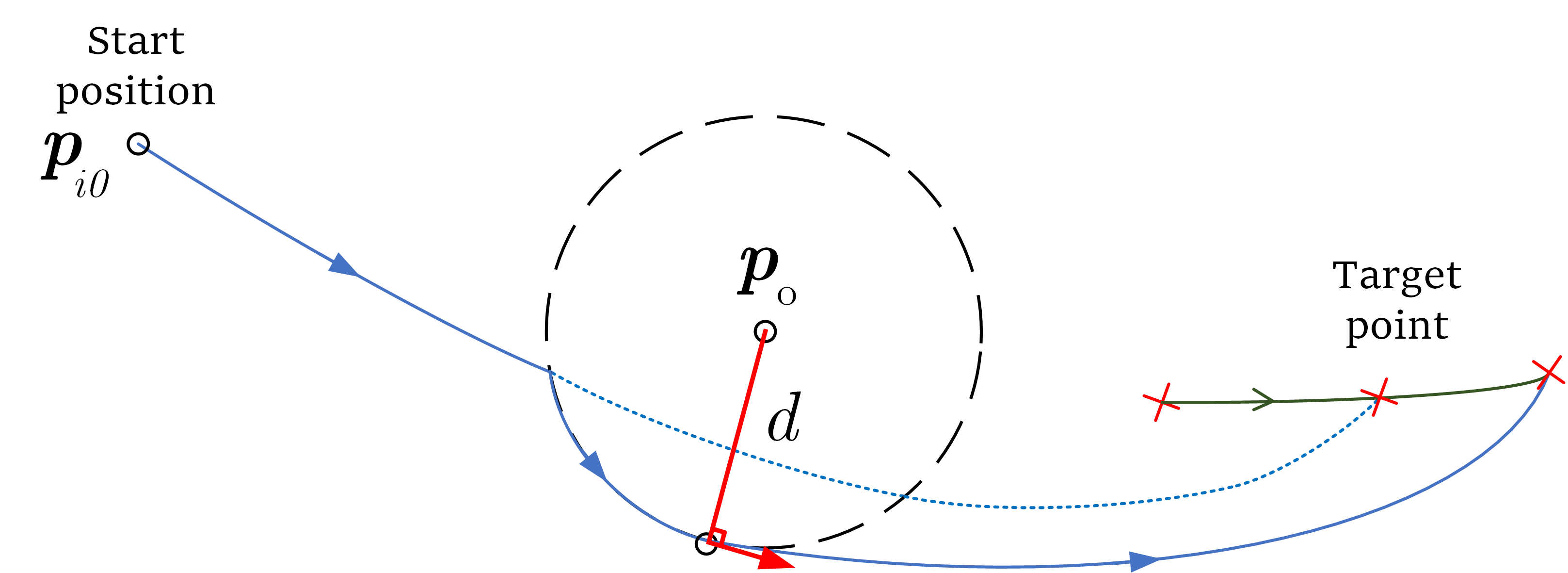} 
\caption{Illustration of the control law (1): the blue dashed-line is the trajectory without the obstacle. The blue line approximates the real trajectory of the agent toward the target point.}
\end{figure*}

\item Case 2: The solution of $\dot{\m{p}}_i = - \Big(\text{sig}\big(\m{p}_i - \hat{\m{p}}_i \big)^\alpha + \beta \text{sign}\big(\m{p}_i - \hat{\m{p}}_i \big) \Big), {\m{p}}_i(t_0) = {\m{p}}_{i0}$ intersects $C$. 
\begin{itemize}
\item Case 2 (a): $\m{p}_i \in C$. In this case, $\zeta = 1$ and the motion of agent $i$ is governed by the control law $\m{u}_i=\big(\m{p}_i - {\m{p}}_o \big) + \m{g}_{io}^\perp$. This control law consists of two terms: the first term steers the agent toward the boundary of $C$ in a finite-time. The second one steers the agent along the orthogonal space of $\m{g}_{io}$. Note that inside $C$, the first control term dominates the second control term (because $k>>1$) and thus it prevents possible collision. The second term changes the bearing vector $\m{g}_{io}$. Thus, the combination of two control laws eventually steers the agent to a position in the boundary of $C$, denoted by $\m{p}_i'$. Two cases may happen here: (i) the solution $S'$ of $\dot{\m{p}}_i = - \Big(\text{sig}\big(\m{p}_i - \hat{\m{p}}_i \big)^\alpha + \beta \text{sign}\big(\m{p}_i - \hat{\m{p}}_i \big) \Big), {\m{p}}_i(0) = \m{p}_i'$ does not intersect the region $C$, then we return to Case 1; (ii) $S'$ intersects the region $C$, then we have Case 2(b). These arguments are illustrated in Figure 1. To show these arguments, consider:
\begin{align*}
\frac{d}{dt}d_{io}^2 &= 2(\m{p}_i-\m{p}_o)^\top \dot{\m{p}_i} \\
&= 2k(\m{p}_i-\m{p}_o)^\top\big(\m{p}_i - {\m{p}}_o \big) + 2(\m{p}_i-\m{p}_o)^\top\m{g}_{io}^\perp\\
&= 2k \|\m{p}_i-\m{p}_o\|^2>0, \forall \m{p}_o \neq \m{p}_i \in C
\end{align*}
This implies $d_{io} \to d$ and $\m{p}_i \to \partial C$ in a finite time. Moreover,
\begin{align*}
\frac{d}{dt} \m{g}_{io} &= -\frac{\m{P}_{\m{g}_{io}}}{d_{io}}\dot{\m{p}_i} = -\frac{\m{P}_{\m{g}_{io}}}{d_{io}}\big(\m{p}_i - {\m{p}}_o \big) - \frac{\m{P}_{\m{g}_{io}}}{d_{io}}\m{g}_{io}^\perp=- \frac{1}{d_{io}} \m{g}_{io}^\perp (\m{g}_{io}^\perp)^\top \m{g}_{io}^\perp = -\frac{1}{d_{io}} \m{g}_{io}^\perp,
\end{align*}
and thus, the bearing vector $\m{g}_{io}$ changes its direction if $\m{p}_i$ enters the set $C$.

\item Case 2 (b) $\m{p}_i \notin C$. Then, the agent $i$ keeps moving toward the target until it enters $C$. After $\m{p}_i \in C$, we analyze as in Case 2(a).
\end{itemize}
\end{itemize}
\end{IEEEproof}

\begin{figure}[h]
\centering
\includegraphics[height=10cm]{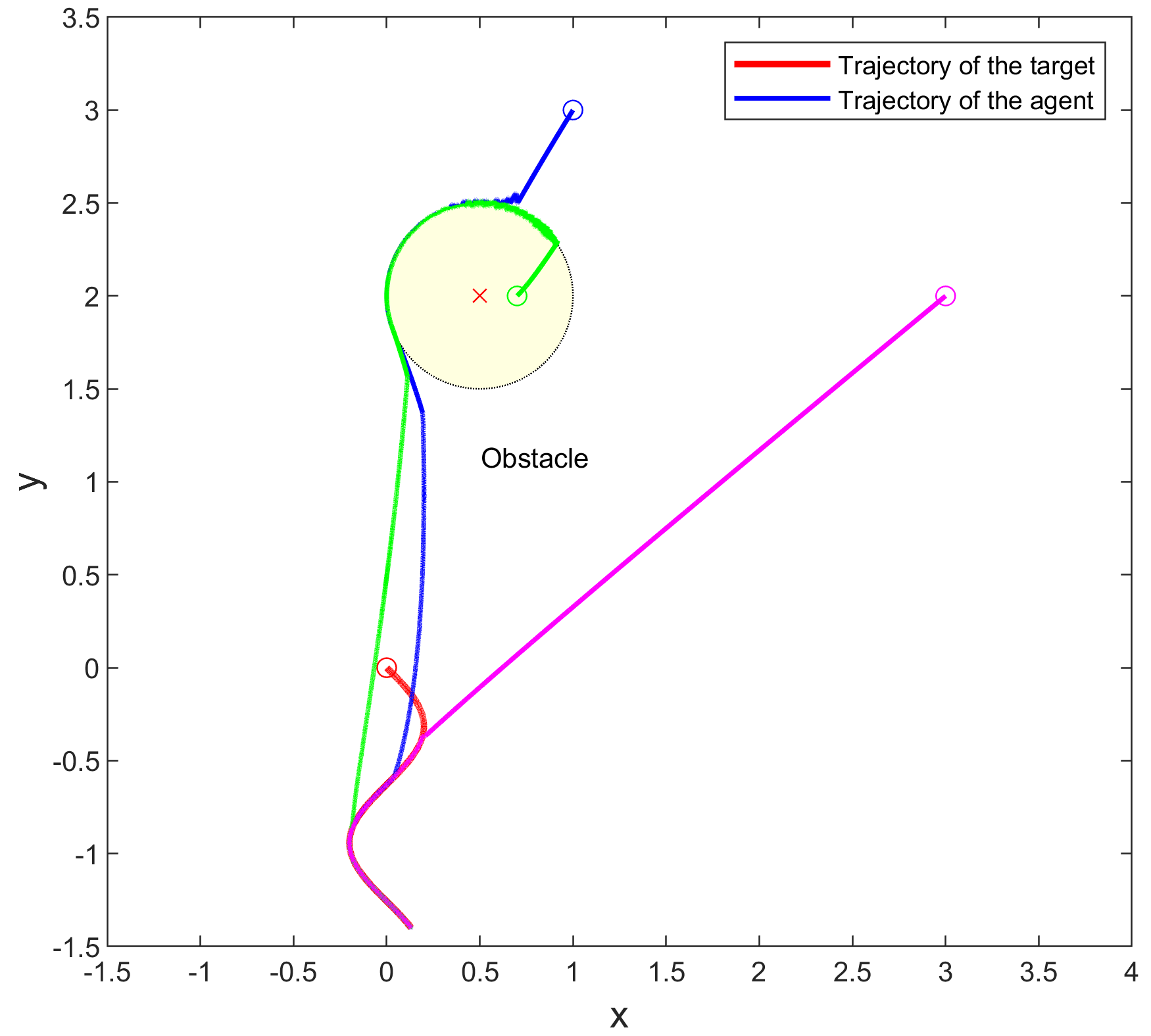}
\caption{Simulation of the control law (1): the agent tracks the target point after finite time.}
\label{fig:1s}
\end{figure}

We simulate the control law \eqref{eq:coav} and demonstrate the result as in Fig.~\ref{fig:1s}. The virtual target point is initially located at $[0,0]^\top$ and it moves with velocity $\hat{\m{v}}_i = [0.2*\cos(t),-0.2]^\top$. The obstacle is located at $[0.5,2]^\top$. The control law's parameters: $d=0.5$ and $k=5$. For three different initial positions of the agent, we have three corresponding trajectories (blue, green, magenta). It can be seen that the agent can track the target point after finite time and achieve collision avoidance.

However, it is noted that the proposed collision scheme will not work if multiple obstacles are presented in the system. 

\begin{Remark} A fixed-time position estimation dynamics can be designed based on \eqref{eq:est1} and \cite{Trinh2020ACC,Polyakov2011nonlinear} as follows:
\begin{align}
\dot{\hat{\m{p}}}_i & = \sum_{j\in \mc{N}_i} \m{P}_{\m{g}_{ij}^*}\Big(\text{sig}(\m{P}_{\m{g}_{ij}^*}(\hat{\m{p}}_j -  \hat{\m{p}}_i))^\alpha + \text{sig}(\m{P}_{\m{g}_{ij}^*}(\hat{\m{p}}_j -  \hat{\m{p}}_i))^\rho\Big) +\gamma_i \sum_{j\in \mc{N}_i} \m{P}_{\m{g}_{ij}^*} \text{sign}\big(\m{P}_{\m{g}_{ij}^*}(\hat{\m{p}}_j - \hat{\m{p}}_i) \big), \label{eq:fixed-time}
\end{align}
where $0<\alpha < 1$ and $\rho > 1$. Note that \eqref{eq:fixed-time} is also applicable for infinitesimally bearing rigid leader-follower formations.
\end{Remark}
\begin{IEEEproof}
Let the leaders be indexed by $1, \ldots, l$ and the followers be indexed by $l+1, \ldots, n$. By noting that the desired moving formation is moving with the same velocity, one may write
\begin{align*}
\dot{\m{p}}^* &= \m{1}_n \otimes \m{v}^*
\end{align*}
where $\m{p}^* = \text{vec}(\m{p}^{l*},\m{p}^{f*})=\text{vec}(\m{p}_1^*, \ldots, \m{p}_n^*)$. Let $\hat{\m{p}}_i$ denote the estimate of the target position of agent $i$, and define $\hat{\m{p}}=\text{vec}(\hat{\m{p}}^l,\hat{\m{p}}^f) = \text{vec}(\hat{\m{p}}_{1}, \ldots, \hat{\m{p}}_n)$. Since the leaders know their positions and has already been at the desired position, $\hat{\m{p}}^l = \m{p}^{l*} = \text{vec}(\m{p}_1^*, \ldots, \m{p}_l^*)$.
Defining the matrix $\m{Z} = \begin{bmatrix}
\m{0}_{dl\times dl} & \m{0}_{dl \times df}\\ \m{0}_{df \times dl} & \m{I}_{df}
\end{bmatrix}$, the equation governing the dynamic of the estimated desired configuration is given as follows:
\begin{align}
\dot{\hat{\m{p}}} = - \m{Z} \bar{\m{H}}^\top \text{diag}(\m{P}_{\m{g}_k^*}) \Big(\text{sig}\big(\text{diag}(\m{P}_{\m{g}_k^*}) \bar{\m{H}} \hat{\m{p}} \big)^\alpha +\text{sig}\big(\text{diag}(\m{P}_{\m{g}_k^*}) \bar{\m{H}} \hat{\m{p}} \big)^\rho +(\m{\Gamma}\otimes \m{I}_d) \text{sign}\big(\text{diag}(\m{P}_{\m{g}_k^*}) \bar{\m{H}} \hat{\m{p}} \big)\Big) + \text{vec}(\m{1}_l \otimes \m{v}^*, \m{0}_{df}), \label{eq:3}
\end{align}
where $\m{\Gamma} = \text{diag}(0, \ldots, 0, \gamma_{l+1} , \ldots, \gamma_n)$. The solution of \eqref{eq:3} is understood in Fillipov sense. 


Without loss of generalization, we will only consider the case $\gamma_i = \gamma_j = \gamma, \forall i = l+1, \ldots, n$ in the following analysis. Consider the Lyapunov function: $V = \frac{1}{2} \|\hat{\m{p}} - \m{p}^*\|^2$, one has $\dot{V}\in^{a.e.} \dot{\tilde{V}}$, and 
\begin{align*}
\dot{\tilde{V}} =& -(\hat{\m{p}} - \m{p}^*)^\top \m{Z} \bar{\m{H}}^\top \text{diag}(\m{P}_{\m{g}_k^*}) \Big(\text{sig}\big(\text{diag}(\m{P}_{\m{g}_k^*}) \bar{\m{H}} \hat{\m{p}} \big)^\alpha +\text{sig}\big(\text{diag}(\m{P}_{\m{g}_k^*}) \bar{\m{H}} \hat{\m{p}} \big)^\rho \\&\qquad +\gamma K[\text{sign}]\big(\text{diag}(\m{P}_{\m{g}_k^*}) \bar{\m{H}}\hat{\m{p}} \big)\Big) -(\hat{\m{p}} - \m{p}^*)^\top[ (\m{1}_l \otimes \m{v}^*, \m{0}_{df}) - \m{1}_n \otimes \m{v}^*] \nonumber\\
=& -\text{vec}(\m{0}_{dl},\hat{\m{p}}^f - \m{p}^{f*})^\top\bar{\m{H}}^\top \text{diag}(\m{P}_{\m{g}_k^*}) \Big(\text{sig}\big(\text{diag}(\m{P}_{\m{g}_k^*}) \bar{\m{H}} \hat{\m{p}} \big)^\alpha +\text{sig}\big(\text{diag}(\m{P}_{\m{g}_k^*}) \bar{\m{H}} \hat{\m{p}} \big)^\rho \\&\qquad +\gamma K[\text{sign}]\big(\text{diag}(\m{P}_{\m{g}_k^*}) \m{H} \hat{\m{p}} \big)\Big) -(\hat{\m{p}}^f - \m{p}^{f*})^\top (\m{1}_f \otimes \m{v}^*) \nonumber\\
=& -(\hat{\m{p}} - \m{p}^*)^\top\bar{\m{H}}^\top \text{diag}(\m{P}_{\m{g}_k^*}) \Big(\text{sig}\big(\text{diag}(\m{P}_{\m{g}_k^*})\bar{\m{H}} \hat{\m{p}} \big)^\alpha +\text{sig}\big(\text{diag}(\m{P}_{\m{g}_k^*}) \bar{\m{H}} \hat{\m{p}} \big)^\rho \\&\qquad +\gamma \text{sign}\big(\text{diag}(\m{P}_{\m{g}_k^*}) \bar{\m{H}} \hat{\m{p}} \big)\Big) -(\hat{\m{p}}^f - \m{p}^{f*})^\top (\m{1}_f \otimes \m{v}^*)\\
=& - |\big(\text{diag}(\m{P}_{\m{g}_k^*})\bar{\m{H}} (\hat{\m{p}}-\m{p}^*)|^{\alpha+1} - |\big(\text{diag}(\m{P}_{\m{g}_k^*}) \bar{\m{H}} (\hat{\m{p}} - \m{p}^*)|^{\rho+1} \\&\qquad - \gamma \|\big(\text{diag}(\m{P}_{\m{g}_k^*}) \bar{\m{H}} (\hat{\m{p}}-\m{p}^*)\|_1 -(\hat{\m{p}}^f - \m{p}^{f*})^\top (\m{1}_f \otimes \m{v}^*)\\
\leq & - \|(\hat{\m{p}}-\m{p}^*)^\top \m{L}_b(\m{p}^*) (\hat{\m{p}}-\m{p}^*)\|^{\frac{\alpha+1}{2}} - \|(\hat{\m{p}}-\m{p}^*)^\top \m{L}_b(\m{p}^*) (\hat{\m{p}}-\m{p}^*)\|^{\frac{\rho+1}{2}} \\&\qquad - \gamma \|\big(\text{diag}(\m{P}_{\m{g}_k^*}) \bar{\m{H}} (\hat{\m{p}}-\m{p}^*)\|_1 + \|\hat{\m{p}}- \m{p}^{*}\| \|\m{v}^*\|
\end{align*}
Note that $\|\big(\text{diag}(\m{P}_{\m{g}_k^*}) \bar{\m{H}} (\hat{\m{p}}-\m{p}^*)\|_1 \geq \|\big(\text{diag}(\m{P}_{\m{g}_k^*}) \bar{\m{H}} (\hat{\m{p}}-\m{p}^*)\| = \Big((\hat{\m{p}}-\m{p}^*)^\top \m{L}_b(\m{p}^*)(\hat{\m{p}}-\m{p}^*)\Big)^{\frac{1}{2}}= (\hat{\m{p}}^f-\m{p}^{f*})^\top \m{L}_{bff}(\m{p}^*)(\hat{\m{p}}^f-\m{p}^{f*})\geq \sqrt{\lambda_{\min}(\m{L}_{bff}(\m{p}^*))}\|\hat{\m{p}}-\m{p}^{*}\| $ .

Thus, if $\gamma$ is chosen such that $\gamma \geq \frac{\|\m{v}^*\|}{ \sqrt{\lambda_{\min}(\m{L}_{bff}(\m{p}^*))}}$, then it is not hard to show that
\begin{align*}
\dot{V} \leq - \chi_1 V^{\frac{\alpha+1}{2}} - \chi_2 V^{\frac{\rho+1}{2}},
\end{align*}
which implies the fixed-time convergence of the estimated configuration to the desired configuration $\m{p}^*$ based on \cite{Polyakov2011nonlinear}.

The fixed-time stability analysis in case of directed acylic leader-follower graphs can be shown by a similar approach as in \cite{Vu2020ICERA,Trinh2020ACC} and will be omitted.
\end{IEEEproof}
\section{Simulation results}
\label{sec:simulation}
\subsection{Simulation 1: Bearing-only control law}
\begin{figure*}
\centering
\subfloat[]{\includegraphics[height = 3.5cm]{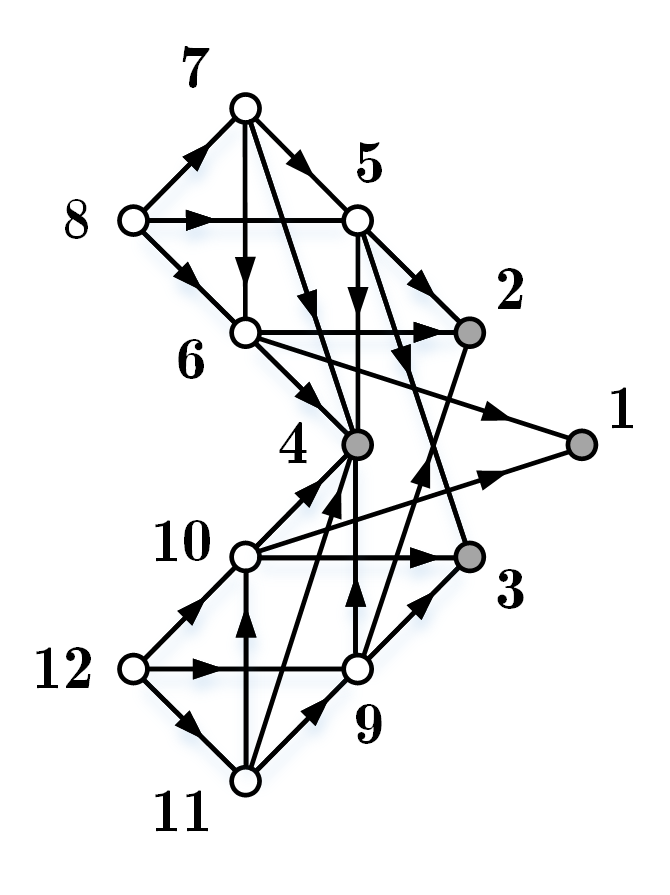}} \qquad
\subfloat[]{\includegraphics[width = 0.785\textwidth]{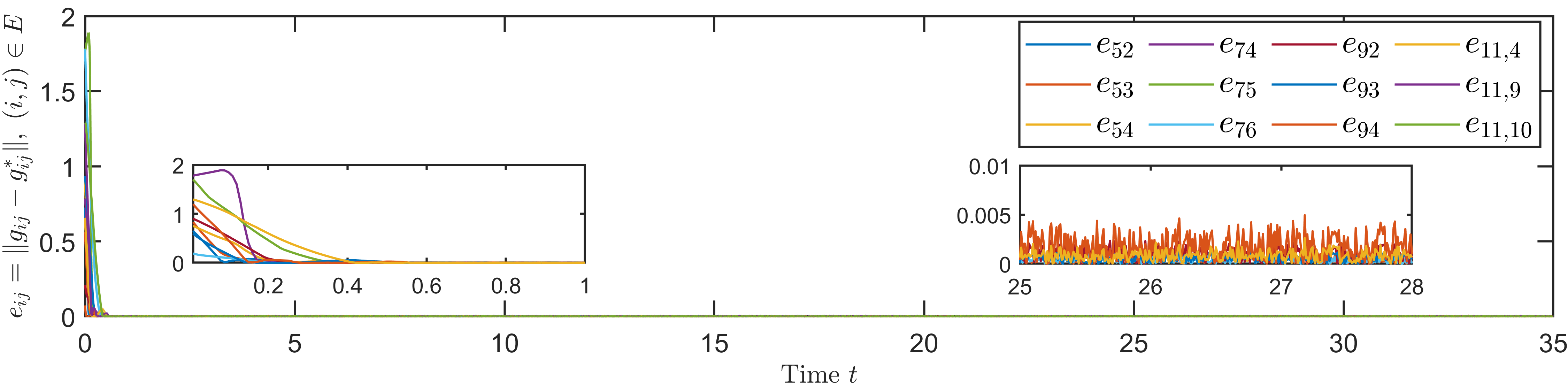}}\\
\subfloat[]{\includegraphics[width = 0.965\textwidth]{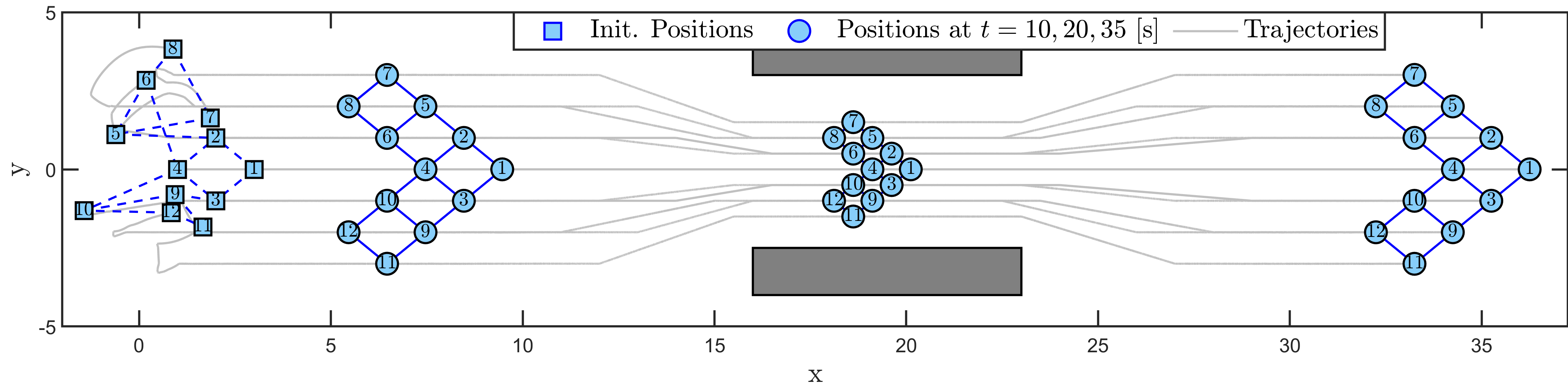}}
\caption{Simulation 1: (a) The acyclic leader-follower graph $\mc{G}$; 
{(b) Bearing errors vs time [s]; (c) Trajectory of the agents.}}
\label{fig:Sim1}
\end{figure*}

Consider a formation of 12 agents with graph as shown in Fig.~\ref{fig:Sim1}(a). Agents 1--4 are leaders, which move with the reference velocity $\m{v}^l$ given as follows
\begin{itemize}
\item For $t \in [0,10]$, $\m{v}^l = \m{1}_4 \otimes \m{f}_1$, where $\m{f}_1 = [1.9-0.14t, 0]^\top$. Leaders move in straight lines along the $x$-axis. At $t=10$, $\m{v}_1(10) = [0.5, 0]^\top$.
\item For $t \in [10,15]$, $\m{v}^l = \m{1}_4 \otimes \m{f}_2 -  \frac{\m{h}}{5\|\m{h}\|}$, where $\m{f}_2 = [0.5, 0]^\top$, $\m{h} = \m{p}^l - \m{1}_4 \otimes \bar{\m{p}}^l$, and $\bar{\m{p}}^l = \frac{1}{4}(\m{p}_1+\m{p}_2+\m{p}_3+\m{p}_4)$ is the geometric center of four leaders. Leaders go along the $x$-axis and downscale the formation's size to fit the alley. At $t=15$s, $\m{v}_1(15) = [0.5, 0]^\top$.
\item For $t \in [15,25]$, $\m{v}^l = \m{1}_4 \otimes \m{f}_3$, where $\m{f}_3 = [0.5+0.05(t-15), 0]^\top$. Leaders move through the alley. At $t=25$s, $\m{v}_1(25) = [1, 0]^\top$.
\item For $t \in [25,30]$, $\m{v}^l = \m{1}_4 \otimes \m{f}_4 + \frac{\m{h}}{5\|\m{h}\|}$, where $\m{f}_4 = [1, 0]^\top$. The formation has passed the alley. Leaders go along the $x$-axis and upscale the formation's size back to normal. At $t=30$s, $\m{v}_1(30) = [1, 0]^\top$.
\item For $t \in [30,35]$, $\m{v}^l = \m{1}_4 \otimes \m{f}_5$, where $\m{f}_5 = [1+0.1(t-30), 0]^\top$. The formation accelerates and continues to move forward along the $x$-axis.
\end{itemize}
The followers adopt the control law \eqref{eq:control-law-1} with $\beta = 2$, $\alpha = 0.5$. The simulation results are depicted in Fig.~\ref{fig:Sim1}(b),(c). The target formation shape is achieved in less than 1 second (see Fig. \ref{fig:Sim1}(b)) and maintained except when the formation rescales its size. Thus, simulation result is consistent with Thm.~\ref{thm:1}.

\subsection{Simulation 2: Target point localization-based control law}
\begin{figure*}
\centering
\subfloat[]{\includegraphics[height = 5.5cm]{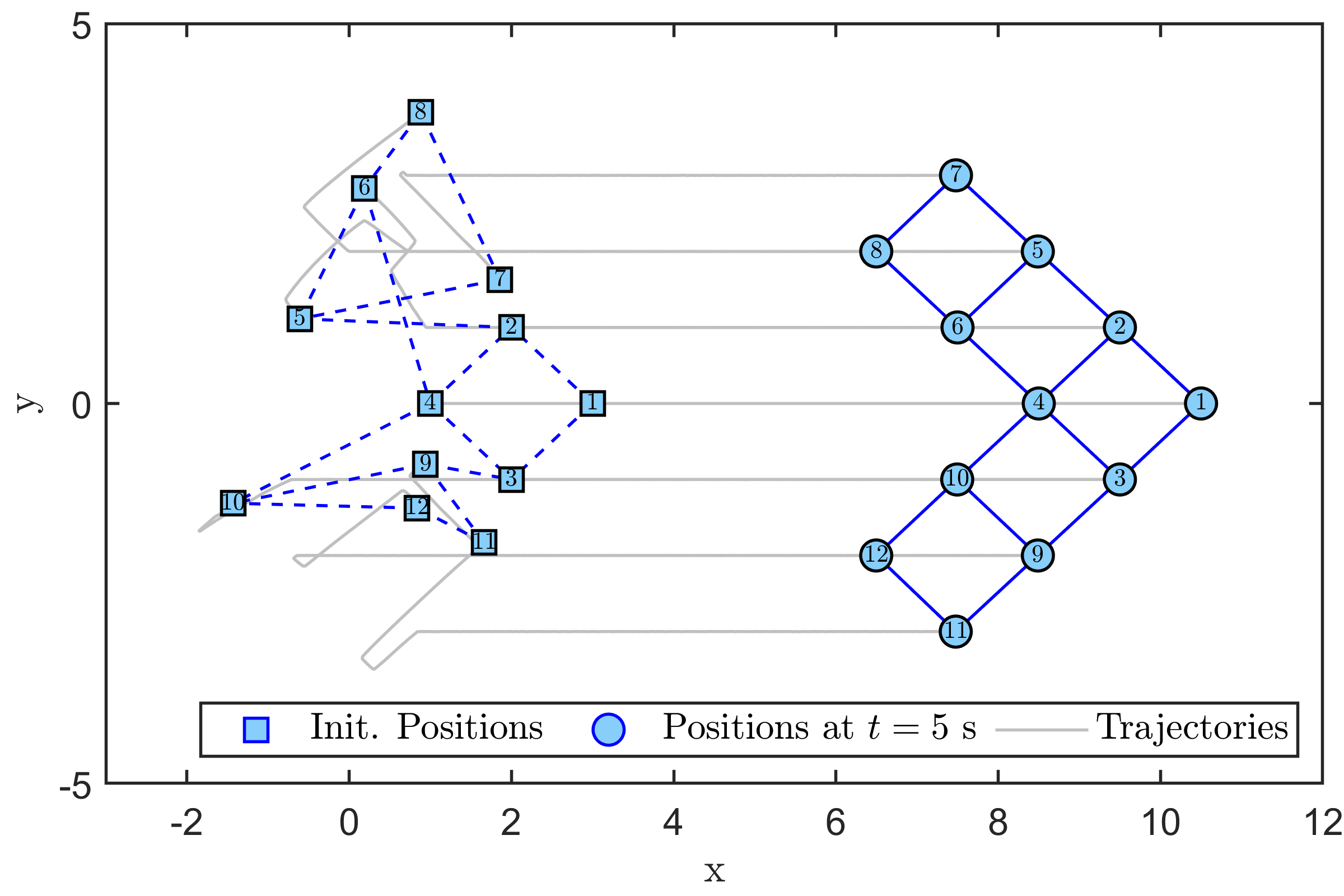}}
\qquad\qquad\quad
\subfloat[]{\includegraphics[height = 5.5cm]{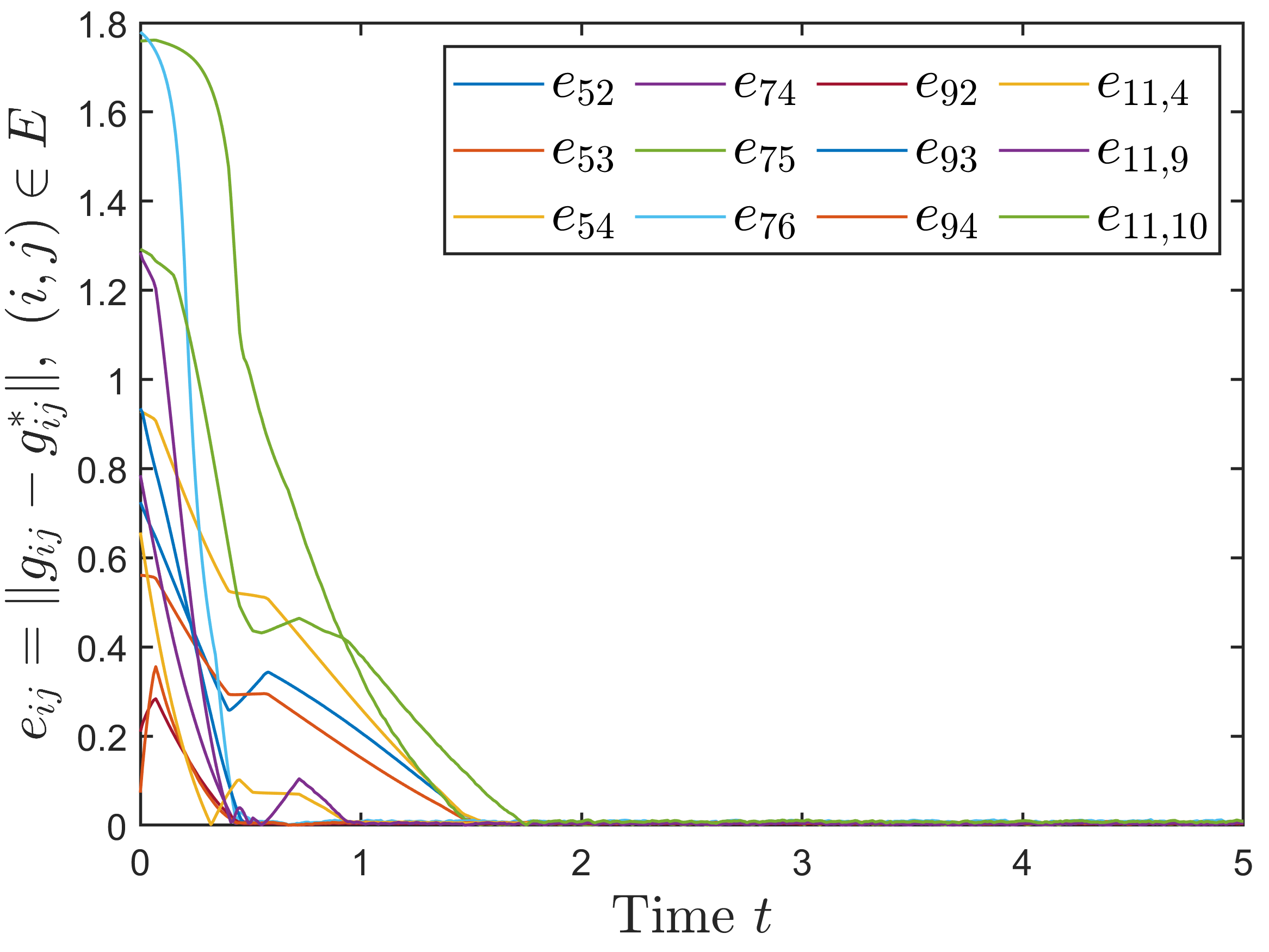}}\qquad

\caption{Simulation 2: (a) Trajectory of the agents. (b) Bearing errors vs time [s].}
\label{fig:Sim2}
\end{figure*}

The same 12-agent formation is simulated under the control law \eqref{eq:est1}--\eqref{eq:est2} for 5 seconds. The initial estimates are randomly selected. The trajectories of the agents are depicted in Fig.~\ref{fig:Sim2}(a). The desired formation is achieved in less than 1 second, which is consistent with Thm.~\ref{thm:est}.
\section{Conclusions}
\label{sec:conclusion}
In this letter, two finite-time bearing-based tracking control laws for acyclic leader-follower formations have been proposed. The analysis partially explains how individuals can follow leaders, who are moving at a time-varying velocity, in collective behaviors such as bird immigration. As suggested in \cite{Nagy2010hierarchical}, it will be interesting to consider the problem with delay and switching in sensing/communication, or when agents can measure only the subtended bearing angles.

\bibliographystyle{IEEEtran}
\bibliography{minh2021}

\end{document}